\documentclass[12pt]{article}
\usepackage{latexsym,amsfonts,amssymb}
\usepackage[dvips]{color}
\usepackage{colordvi,multicol}
\usepackage{amsmath}
\usepackage{graphicx}
\usepackage{pdflscape}
\usepackage{rotating}
\usepackage{breqn}

\newtheorem{thm}{Theorem}[section]

\newtheorem{lem}[thm]{Lemma}
\newtheorem{pro}[thm]{Proposition}

\newtheorem{rem}[thm]{Remark}

\newtheorem{con}[thm]{Conjecture}

\date{}
\begin{document}

\title{\bf Using Sums-of-Squares to Prove Gaussian Product
Inequalities}
 \author{Oliver Russell and Wei Sun\\ \\ \\
  {\small Department of Mathematics and Statistics}\\
    {\small Concordia University, Canada}\\ \\
{\small o\_russel@live.concordia.ca,\ \ \ \ wei.sun@concordia.ca}}

\maketitle

\begin{abstract}

\noindent The long-standing Gaussian product inequality (GPI) conjecture states that
$E [\prod_{j=1}^{n}X_j^{2m_j}]\geq\prod_{j=1}^{n}E[X_j^{2m_j}]$ for any centered Gaussian random vector $(X_1,\dots,X_n)$ and $m_1,\dots,m_n\in\mathbb{N}$. In this paper, we describe a computational algorithm involving sums-of-squares representations of multivariate polynomials that can be used to resolve the GPI conjecture. To exhibit the power of this novel method, we apply it to prove two new GPIs: $E[X_1^{2m_1}X_2^{6}X_3^{4}]\ge E[X_1^{2m_1}]E[X_2^{6}]E[X_3^{4}]$ and $E[X_1^{2m_1}X_2^{2}X_3^{2}X_4^{2}]\ge E[X_1^{2m_1}]E[X_2^{2}]E[X_3^{2}]E[X_4^{2}]$.
\end{abstract}



\section{Introduction}

Gaussian distributions play a fundamental role in both probability and statistics (see \cite{T} and references therein). Over the past decade, there has been a surge of interest in the study of inequalities related to Gaussian distributions. This was sparked by Royen's famous proof of the long-conjectured Gaussian correlation inequality (GCI) \cite{Royen,LM}. Gaussian inequalities are valuable tools to the study of small ball probabilities (see, for example, Li \cite{Li} and Shao \cite{Shao}) and the zeros of random polynomials (see, for example, Li and Shao \cite{LS}), among others. In this paper, we will investigate the Gaussian product inequality (GPI), another long-standing conjecture and one of the late Wenbo Li's open problems (see Shao \cite{Shao2}).

In its most general form, the GPI conjecture (see Li and Wei \cite{LW12}) states that for any non-negative real numbers $\alpha_j$, $j=1,\ldots,{n}$, and any ${n}$-dimensional real-valued centered Gaussian random vector $(X_1,\dots,X_{n})$,
\begin{eqnarray}\label{LW-inequ}
E \left[\prod_{j=1}^{n}|X_j|^{\alpha_j}\right]\geq \prod_{j=1}^{n}E[|X_j|^{\alpha_j}].
\end{eqnarray}
Verification of this inequality has immediate and far-reaching consequences. For example (see Malicet et al. \cite{MNPP16}), if (\ref{LW-inequ}) holds when the $\alpha_j$'s are any equal positive  even integers, then the `real linear polarization constant' conjecture raised by Ben\'{\i}tem  et al. \cite{BST98} is true, and (\ref{LW-inequ}) is deeply linked to the celebrated $U$-conjecture by Kagan, Linnik and Rao \cite{KLR1973}.

The GPI conjecture has proven very challenging to solve, but not for lack of trying. Some progress has been made in the form of partial results. In \cite{Fr08}, Frenkel used algebraic methods to prove (\ref{LW-inequ}) for the case that $\alpha_j=2$ for each $j=1,\ldots,{n}$. In \cite[Theorem 3.2]{LHS}, Lan et al. used the Gaussian hypergeometric functions to prove the following GPI: for any $m_1,m_2\in\mathbb{N}$ and any centered Gaussian random vector $(X_1,X_2,X_3)$,
$$
E[X_1^{2m_1}X_2^{2m_2}X_3^{2m_2}]\ge E[X_1^{2m_1}]E[X_2^{2m_2}]E[X_3^{2m_2}].
$$
Lan et al.'s proof of the $3$-dimensional GPI brought renewed hope to solving this difficult problem.

Throughout this paper, any Gaussian random variable is assumed to be real-valued and non-degenerate, i.e., has positive variance. We will study the following conjecture.
\begin{con}\label{con1} Let $n\ge 3$ and $m_1,\dots,m_n\in\mathbb{N}$.
For any centered Gaussian random vector $(X_1,\dots,X_n)$,
\begin{equation}\label{14AA}
E \left[\prod_{j=1}^{n}X_j^{2m_j}\right]\geq\prod_{j=1}^{n}E[X_j^{2m_j}].
\end{equation}
The equality holds if and only if $X_1,\dots,X_n$ are independent.
\end{con}

Most recently, using combinatorial methods we proved Conjecture \ref{con1} for $n=3$ when $m_1=1,m_2=2,m_3\in\mathbb{N}$ and $m_1=1,m_2=3,m_3\in\mathbb{N}$ (see \cite[Theorem 3.5]{RusSun}). In addition, in the same paper, we proved the GPI (for any $n\in\mathbb{N}$) when all correlations are non-negative.

\begin{pro}\label{pro1}(\cite[Lemma 2.3]{RusSun})
Let $(X_1,\dots,X_n)$ be a centered Gaussian random vector such that $E[X_iX_j]\ge 0$ for any $i\not=j$. Then, (\ref{14AA}) holds and
\begin{eqnarray*}
E \left[\prod_{j=1}^{n}X_j^{2m_j}\right]\geq
E \left[\prod_{j=1}^kX_j^{2m_j}\right]E\left[\prod_{j=k+1}^{n}X_j^{2m_j}\right],\ \ \ \ \forall 1\le k\le n-1.
\end{eqnarray*}

\end{pro}

\begin{rem} (i) Genest and Ouimet \cite{GenOui} proved that (\ref{14AA}) holds under the assumption that there exists a matrix $C\in[0,\infty)^{n\times n}$ such that $(X_1,\dots,X_n)=(U_1,\dots,U_n)C$ in law, where $(U_1,\dots,U_n)$
is an $n$-dimensional standard Gaussian random vector. Recently, Edelmann et al. \cite{Edel} used a different method to extend \cite[Lemma 2.3]{RusSun} to the multivariate gamma distribution. It is interesting to point out that Edelmann et al. [3] also demonstrated a strong
connection between the GCI and the GPI.

\noindent (ii) Note that the product inequality (\ref{14AA}) may not hold for non-Gaussian random vectors. Let $V_1$ and $V_2$ be independent Rademacher random variables, i.e., $P(V_i=1)=P(V_i=-1)=\frac{1}{2}$, $i=1,2$. Define $Y_1=V_1+0.9V_2$ and $Y_2=V_1-0.9V_2$. Then, $(Y_1,Y_2)$ is a centered random vector with
$$
{\rm Cov}(Y_1,Y_2)=E[V_1^2-0.81V_2^2]=0.19.
$$
We have
$$
E[Y_1^2Y_2^2]=E[(V_1^2-0.81V_2^2)^2]=(0.19)^2,
$$
and
$$
E[Y_1^2]E[Y_2^2]=(1.81)^2.
$$
Hence
$$
E[Y_1^2Y_2^2]<E[Y_1^2]E[Y_2^2].
$$

\end{rem}

In this paper, we develop an efficient computational algorithm that produces \textit{exact} sums-of-squares (SOS) polynomials to tackle the GPI. We describe this method in Section 2 and use it to rigorously prove two previously unsolved special cases of the GPI conjecture with fixed exponents in Section 3. Then, in Section 4, we reveal the true power of the SOS method by extending these special cases to the stronger result where one exponent is unbounded. Fully solving the GPI conjecture is most difficult when correlations can be negative. In theory, our algorithm is applicable to \textit{any} GPI of the form (\ref{14AA}) with dimension and all but one exponent fixed, and is therefore the first universal method for solving the GPI.

\section{The SOS method of solving the GPI}\setcounter{equation}{0}

An SOS representation of a polynomial is of the form $\sum_{i=1}^p f_i^2$, where  the $f_i$'s are real-coefficient polynomials. It is clear that any polynomial with an SOS representation is necessarily non-negative.

\begin{lem} Let $(X_1,\dots,X_n)$ be a centered Gaussian random vector.
Denote by $\Lambda$ the covariance matrix of $(X_1,\dots,X_n)$ and $c_{m_1,\dots,m_n}$ the coefficient corresponding to the term $t_1^{2 m_1}\cdots t_n^{2 m_n}$ of the polynomial
$$
G(t_1,\dots,t_n)=\left(\sum_{k,l=1}^n\Lambda_{kl}t_kt_l\right)^{\sum_{j=1}^nm_j},\ \ \ \ t_1,\dots,t_n\in \mathbb{R}.
$$
Then,
\begin{equation}\label{14B}
E \left[\prod_{j=1}^{n}X_j^{2m_j}\right]=\frac{\prod_{j=1}^n(2m_j)!}{2^{\sum_{j=1}^nm_j}(\sum_{j=1}^nm_j)!}\cdot c_{m_1,\dots,m_n}.
\end{equation}
\end{lem}
\noindent{\bf Proof.}\ \
We use $'$ to denote the transpose of a matrix or a vector. Define
$$
f_{\Lambda}(x)=\frac{1}{\sqrt{(2\pi)^n{\rm det}(\Lambda)}}e^{-\frac{1}{2}x'{\Lambda}^{-1}x},\ \ x=(x_1,\dots,x_n)\in \mathbb{R}^n,
$$
and
$$
\Psi_{\Lambda}(t)= \int_{\mathbb{R}^n}e^{it'x}f_{\Lambda}(x)dx,\ \ \ \ t=(t_1,\dots,t_n)\in \mathbb{R}^n.
$$
Then,
$$
 \Psi_{\Lambda}(t)=\exp\left(-\frac{1}{2}t'\Lambda t\right).
$$

For $\beta=(\beta_1,\dots,\beta_n)$ with $\beta_j\in\mathbb{N}\cup\{0\}$, we have
$$
\partial^{\beta}_t\Psi_{\Lambda}(0)=i^{\sum_{j=1}^n\beta_j}\int_{\mathbb{R}^n}\left\{\prod_{j=1}^nx_j^{\beta_j}\right\}f_{\Lambda}(x)dx.
$$
Hence,
\begin{eqnarray*}
E \left[\prod_{j=1}^{n}X_j^{2m_j}\right]&=&
(-1)^{\sum_{j=1}^nm_j}\partial^{(2m_1,\dots,2m_n)}_t\exp\left(-\frac{1}{2}t'\Lambda t\right)(0)\nonumber\\
&=&\frac{1}{2^{\sum_{j=1}^nm_j}(\sum_{j=1}^nm_j)!}\cdot
\partial^{(2m_1,\dots,2m_n)}_t\left(\sum_{k,l=1}^n\Lambda_{kl}t_kt_l\right)^{\sum_{j=1}^nm_j}(0).
\end{eqnarray*}
Therefore, the proof is complete. \hfill \fbox

Let $U_j$, $1\le j\le n$, be independent standard Gaussian random variables. Define
$$
X_k=\sum_{j=1}^{n}x_{kj}U_j,\ \ 1\le k\le n,
$$
where each $x_{kj}\in\mathbb{R}$, $1\leq k,j\leq n$. Then, we have
\begin{eqnarray}\label{19A}
&&\Lambda_{kk}=\sum_{j=1}^{n}x^2_{kj},\ \ 1\le k\le n,\nonumber\\
&&\Lambda_{kl}=\sum_{j=1}^{n}x_{kj}x_{lj},\ \ 1\le k< l\le n.
\end{eqnarray}
 Define
\begin{eqnarray*}
F_{m_1,\dots,m_n}(\Lambda)&=&E \left[\prod_{j=1}^{n}X_j^{2m_j}\right]-\prod_{j=1}^{n}E[X_j^{2m_j}]\\
&=&E \left[\prod_{j=1}^{n}X_j^{2m_j}\right]-\prod_{k=1}^n[(2m_k-1)!!\Lambda_{kk}^{m_k}].
\end{eqnarray*}
 By (\ref{14B}) and (\ref{19A}), it is easy to see that $F_{m_1,\dots,m_n}(\Lambda)$ can be expressed as a polynomial of the $x_{ij}$'s, say $F_{m_1,\dots,m_n}(x_{11},\dots,x_{1n},\dots,x_{n1},\dots,x_{nn})$.

By (\ref{14B}) and using the \textbf{Mathematica} functions \texttt{Expand[]} and \texttt{Coefficient[]}, we get
\begin{eqnarray}\label{19B}
{\rm Poly}&:=&{\rm \texttt{Expand}}\left[\left(\sum_{k,l=1}^n\Lambda_{kl}t_kt_l\right)^{\sum_{j=1}^nm_j}\right],\nonumber\\
F_{m_1,\dots,m_n}(\Lambda)
&=&{\rm \texttt{Expand}}\left[\frac{\prod_{j=1}^n(2m_j)!}{2^{\sum_{j=1}^nm_j}(\sum_{j=1}^nm_j)!}\cdot
   {\rm \texttt{Coefficient}}\left[{\rm Poly},
    \prod_{j=1}^nt_j^{2 m_j}\right] \right.\nonumber\\
&&\ \ \ \ \ \ \ \ \ \ \ \ \ \ \left.- \prod_{k=1}^n\left[(2m_k-1)!!\Lambda_{kk}^{m_k}\right]\right].
\end{eqnarray}
Combining (\ref{19A}) and (\ref{19B}), we have an algorithm to get the expression of \\ $F_{m_1,\dots,m_n}(x_{11},\dots,x_{1n},\dots,x_{n1},\dots,x_{nn})$.
Suppose we consider (\ref{19B}) for fixed $n, m_1,\dots,m_n$. Applying the \textbf{Macaulay2} package \texttt{`SumsOfSquares'} \cite{SOS} to this polynomial, we may attempt to obtain an SOS decomposition $\sum_{i=1}^p c_i f_i^2$  of $F_{m_1,\dots,m_n}$, where the $c_i$'s are positive rational numbers and  the $f_i$'s are rational-coefficient polynomials. Although not every rational-coefficient SOS polynomial necessarily has such a decomposition (see Scheiderer \cite[Theorem 2.1]{Sch16}), this software aims to produce one regardless (see Peyrl and Parrilo \cite{PP}). If it is obtained, then (\ref{14AA}) is verified for this case.

This package is very user-friendly and, along with the semi-definite programming package it uses, comes pre-installed in the newest versions of \textbf{Macaulay2}. By using a ``mixed symbolic-numerical approach" \cite{PP}, \texttt{`SumsOfSquares'} takes advantage of the speed of approximate numerical calculations, yet still produces a final SOS decomposition that is exact (\textit{not} an approximation). This SOS polynomial can then be expanded and checked to match the original $F_{m_1,\dots,m_n}$ using the \texttt{value()} function.

To further increase the efficiency of our method, we will use a trick to reduce the degree of the polynomial $F_{m_1,\dots,m_n}$. In this way, calculation of the SOS decomposition becomes significantly less computationally intensive and much faster. This technique resides in a key result proved in our previous paper which we restate here for convenience.

\begin{lem}\label{lem2}(\cite[Lemma 2.1]{RusSun}) Let $n\ge 3$ and $m_1,\dots,m_n\in \mathbb{N}$. If for any centered Gaussian random vector $(Y_1,\dots,Y_n)$ with $Y_n=\alpha_1Y_1+\cdots+\alpha_{n-1}Y_{n-1}$ for some constants
$\alpha_1,\dots,\alpha_{n-1}$,
\begin{equation}\label{eqnYn2kll}
    E \left[\left\{\prod_{j=1}^{n-1}Y_j^{2m_j}\right\}Y_n^{2k}\right]\ge \left\{\prod_{j=1}^{n-1}E[Y_j^{2m_j}]\right\}E[Y_n^{2k}],\ \ \ \ 0\leq k\leq m_n,
\end{equation}
then for any centered Gaussian random vector $(X_1,\dots,X_n)$,
\begin{equation}\label{14AAll}
    E \left[\prod_{j=1}^{n}X_j^{2m_j}\right]\ge \prod_{j=1}^{n}E[X_j^{2m_j}].
\end{equation}
Additionally, if inequality (\ref{eqnYn2kll}) is strict when $k=m_n$, then the equality sign of (\ref{14AAll}) holds only if $X_n$ is independent of $X_1,\dots,X_{n-1}$.

\end{lem}

\section{Applications of the SOS method}\setcounter{equation}{0}

First, we will verify (\ref{14AA}) for the case $n=3, m_1=4, m_2=3, m_3=2$.

\begin{thm}\label{thm234}
For any centered Gaussian random vector $(X_1,X_2,X_3)$,
$$
E[X_1^{8} X_2^{6} X_3^{4}]\ge E[X_1^{8}] E[X_2^{6}] E[X_3^{4}].
$$
The equality holds if and only if $X_1,X_2,X_3$ are independent.
\end{thm}
\noindent{\bf Proof.}\ \
Let $X_3$ be a linear combination of $X_1$ and $X_2$, and let $U_1$ and $U_2$ be independent standard Gaussian random variables. Then, without loss of generality, we may write
$$
X_1=U_1, \,\, X_2=aU_1+U_2, \,\, X_3=bU_1+U_2, \,\,\, a,b \in \mathbb{R}.
$$
By Lemma \ref{lem2}, we need only show that $F_{4,3,1}$ and $F_{4,3,2}$ are strictly positive. By \cite[Theorem 3.5]{RusSun}, $F_{4,3,1}>0$. We will complete the proof by giving an SOS decomposition of $F_{4,3,2}$. By (\ref{19B}), we get
\begin{dmath*}
F_{4,3,2}(a,b)=94500 + 1474200 \,a^2 + 2324700 \,a^4 + 400680 \,a^6 + 2381400 \,a \,b + 12474000 \,a^3 \,b + 9729720 \,a^5 \,b + 585900 \,b^2 + 14004900 \,a^2 \,b^2 + 36458100 \,a^4 \,b^2 + 12152700 \,a^6 \,b^2 + 3742200 \,a \,b^3 + 32432400 \,a^3 \,b^3 + 48648600 \,a^5 \,b^3 + 151200 \,b^4 + 6066900 \,a^2 \,b^4 + 30391200 \,a^4 \,b^4 + 34454700 \,a^6 \,b^4
\end{dmath*}
\begin{scriptsize}\begin{dmath*}
=\frac{145216503}{4} \bigg(\frac{60146041}{290433006}\,a^{3}b+a^{2}b^{2}+\frac{2953715}{145216503}\,a^{2}+\frac{88269419}{290433006}\,a\,b-\frac{244387}{12627522}\,b^{2}+\frac{54414145}{2468680551} \bigg)^2 + 34454700\bigg(a^{3}b^{2}-\frac{8427157}{275637600}\,a^{3}+\frac{134448359}{275637600}\,a^{2}b-\frac{875989}{10208800}\,a\,b^{2}+\frac{639317}{9844200}\,a-\frac{97603}{5512752}\,b \bigg)^2 + \frac{43421589679843919}{2205100800} \bigg(\frac{2207851366454963}{43421589679843919}\,a^{3}+a^{2}b+\frac{16239747161830977}{43421589679843919}\,a\,b^{2}+\frac{9538031563407116}{43421589679843919}\,a+\frac{55963631836960250}{738167024557346623}\,b \bigg)^2 + \frac{29513864367243803}{2323464048} \bigg(a^{3}b+\frac{9815490958053974}{29513864367243803}\,a^{2}-\frac{2663681330163887}{29513864367243803}\,a\,b+\frac{36010094112641}{29513864367243803}\,b^{2}-\frac{15346278961059626}{501735694243144651} \bigg)^2 + \frac{775345651901116499475901}{173686358719375676} \bigg(-\frac{51214758446985696098233}{775345651901116499475901}\,a^{3}+a\,b^{2}-\frac{124746121201440746832731}{13180876082318980491090317}\,a+\frac{6413514425190010115326451}{26361752164637960982180634}\,b \bigg)^2 + \frac{3557651776301827092683371}{1003471388486289302} \bigg(\frac{5097266129812145832080311}{14230607105207308370733484}\,a^{2}+a\,b+\frac{1562593143932213174576929}{14230607105207308370733484}\,b^{2}+\frac{1355975034048165669823823}{14230607105207308370733484} \bigg)^2 + \frac{69744372474088890758457139816341}{113844856841658466965867872} \bigg(a^{2}-\frac{1715759360111955337524806751893}{131739370228834571432641264097533}\,b^{2}-\frac{186306503137515128525276484100099}{1185654332059511142893771376877797} \bigg)^2 + \frac{497216566987657482387717235351591}{896299573597690673394141556} \bigg(\frac{265821781739646648845732660243440}{497216566987657482387717235351591}\,b-\frac{211256870027940529550680035362988}{497216566987657482387717235351591}\,a^{3}+a \bigg)^2 + \frac{98420585988442784935014259601772304508}{497216566987657482387717235351591} \bigg(a^{3}+\frac{2234800582191780362981720603629611853225}{13385199694428218751161939305841033413088}\,b \bigg)^2 + \frac{212159700340573881370359421737581680593}{2239569293890187714354901489658061} \bigg(b^{2}-\frac{40248036588293144056679077001926201625}{121234114480327932211633955278618103196} \bigg)^2 + \frac{94193588538517108050774757159729177812091894697}{1820387158442237750158023745594380544179968} \bigg( b \bigg)^2 + \frac{1012078525061381023480077400990087498959489065}{140146636339259089636648852302082527294576} \bigg( 1 \bigg)^2\\
>0,
\end{dmath*}\end{scriptsize}

\noindent where the second expression (i.e. the SOS decomposition) is obtained by an application of \texttt{`SumsOfSquares'} to the first expression. \hfill \fbox

Next, we will verify (\ref{14AA}) for the case $n=4, m_1=2, m_2=1, m_3=1, m_4=1$.

\begin{thm}\label{thm1112}
For any centered Gaussian random vector $(X_1,X_2,X_3,X_4)$,
\begin{equation}\label{ineq1112}
E[X_1^{4} X_2^{2} X_3^{2} X_4^{2}]\ge E[X_1^{4}] E[X_2^{2}] E[X_3^{2}] E[X_4^{2}].
\end{equation}
The equality holds if and only if $X_1,X_2,X_3,X_4$ are independent.
\end{thm}
\noindent{\bf Proof.}\ \
Let $X_4$ be a linear combination of $X_1$, $X_2$ and $X_3$, and let $U_1$, $U_2$ and $U_3$ be independent standard Gaussian random variables. Then, without loss of generality, we may consider 3 constructions of $X_1,\dots,X_4$ in terms of $U_1,U_2,U_3$:

\noindent \textit{Case 1}:
$$
X_1=U_1, \,\, X_2=aU_1+U_2, \,\, X_3=bU_1+cU_2+U_3, \,\, X_4=dU_1+eU_2+U_3, \,\,\, a,b,c,d,e \in \mathbb{R}.
$$
\textit{Case 2}:
$$
X_1=U_1, \,\, X_2=aU_1+U_2, \,\, X_3=bU_1+U_2, \,\, X_4=cU_1+dU_2+U_3, \,\,\, a,b,c,d \in \mathbb{R}.
$$
\textit{Case 3}:
$$
X_1=U_1, \,\, X_2=aU_1+U_2, \,\, X_3=bU_1+U_2, \,\, X_4=cU_1+U_2, \,\,\, a,b,c \in \mathbb{R}.
$$
By Lemma \ref{lem2}, we need only show that, in each of these $3$ cases, $F_{2,1,1,1}>0$.

\noindent \textit{Case 1}:
\begin{dmath*}
F_{2,1,1,1}(a,b,c,d,e)=6 + 42 \,a^2 + 12 \,b^2 + 102 \,a^2 \,b^2 + 60 \,a \,b \,c + 6 \,c^2 + 12 \,a^2 \,c^2 + 60 \,b \,d + 420 \,a^2 \,b \,d + 120 \,a \,c \,d + 12 \,d^2 + 102 \,a^2 \,d^2 + 102 \,b^2 \,d^2 + 942 \,a^2 \,b^2 \,d^2 + 420 \,a \,b \,c \,d^2 + 42 \,c^2 \,d^2 + 102 \,a^2 \,c^2 \,d^2 + 120 \,a \,b \,e + 36 \,c \,e + 60 \,a^2 \,c \,e + 60 \,a \,d \,e + 420 \,a \,b^2 \,d \,e + 180 \,b \,c \,d \,e + 420 \,a^2 \,b \,c \,d \,e + 180 \,a \,c^2 \,d \,e + 6 \,e^2 + 12 \,a^2 \,e^2 + 42 \,b^2 \,e^2 + 102 \,a^2 \,b^2 \,e^2 + 180 \,a \,b \,c \,e^2 + 42 \,c^2 \,e^2 + 42 \,a^2 \,c^2 \,e^2
\end{dmath*}
\begin{footnotesize}\begin{dmath*}
=942 \bigg( a\,b\,d+\frac{1169}{7536}\,a\,c\,e+\frac{1159}{7536}\,c\,d+\frac{1159}{7536}\,b\,e+\frac{1217}{7536}\,a \bigg)^2 + 102 \bigg( a\,b+\frac{463}{816}\,a\,d+\frac{29}{136}\,c+\frac{7}{51}\,e \bigg)^2 + 102 \bigg( \frac{511}{816}\,a\,c\,d+a\,b\,e+\frac{521}{816}\,b\,d+\frac{71}{136}\,c\,e+\frac{11}{68} \bigg)^2 + \frac{451487}{6528} \bigg( a\,d+\frac{10830}{451487}\,c+\frac{90128}{451487}\,e \bigg)^2 + \frac{404735}{6528} \bigg( a\,c\,d+\frac{521}{1327}\,b\,d+\frac{426}{1327}\,c\,e+\frac{132}{1327} \bigg)^2 + \frac{269975}{5308} \bigg( b\,d+\frac{25609}{107990}\,c\,e+\frac{10397}{53995} \bigg)^2 + \frac{1188815}{60288} \bigg( \frac{860713}{1188815}\,a\,c\,e+\frac{834623}{1188815}\,c\,d+b\,e+\frac{134377}{1188815}\,a \bigg)^2 + \frac{81700883}{4755260} \bigg( \frac{17007629}{163401766}\,a\,c\,e+\frac{6315719}{163401766}\,c\,d+a \bigg)^2 + 12 \bigg( \frac{1}{4}\,a\,c+a\,e+\frac{31}{96}\,b+\frac{11}{16}\,d \bigg)^2 + \frac{45}{4} \bigg( a\,c+\frac{233}{360}\,b+\frac{29}{180}\,d \bigg)^2 + \frac{13038091541}{1307214128} \bigg( \frac{116381935}{277406203}\,a\,c\,e+c\,d \bigg)^2 + \frac{7828043977}{1109624812} \bigg( a\,c\,e \bigg)^2 + \frac{2173}{360} \bigg( b+\frac{9829}{17384}\,d \bigg)^2 + \frac{855377}{172784} \bigg( c\,e+\frac{230438}{855377} \bigg)^2 + \frac{13706281}{3337728} \bigg( d \bigg)^2 + \frac{597018}{451487} \bigg( c+\frac{782927}{2388072}\,e \bigg)^2 + \frac{11273665}{9552288} \bigg( e \bigg)^2 + \frac{403986}{855377} \bigg( 1 \bigg)^2\\
>0,
\end{dmath*}\end{footnotesize}

\noindent \textit{Case 2}:
\begin{dmath*}
F_{2,1,1,1}(a,b,c,d)=6 + 12 \,a^2 + 60 \,a \,b + 12 \,b^2 + 102 \,a^2 \,b^2 + 42 \,c^2 + 102 \,a^2 \,c^2 + 420 \,a \,b \,c^2 + 102 \,b^2 \,c^2 + 942 \,a^2 \,b^2 \,c^2 + 180 \,a \,c \,d + 180 \,b \,c \,d + 420 \,a^2 \,b \,c \,d + 420 \,a \,b^2 \,c \,d+ 42 \,d^2 + 42 \,a^2 \,d^2 + 180 \,a \,b \,d^2 + 42 \,b^2 \,d^2 + 102 \,a^2 \,b^2 \,d^2
\end{dmath*}
\begin{dmath*}
=942 \bigg( a\,b\,c+\frac{583}{3768}\,a\,d+\frac{583}{3768}\,b\,d+\frac{583}{3768}\,c \bigg)^2 + 102 \bigg( a\,b+\frac{43}{204} \bigg)^2 + 102 \bigg( a\,b\,d+\frac{257}{408}\,a\,c+\frac{257}{408}\,b\,c+\frac{25}{48}\,d \bigg)^2 + \frac{100415}{1632} \bigg( a\,c+\frac{257}{665}\,b\,c+\frac{85}{266}\,d \bigg)^2 + \frac{69611}{1330} \bigg( b\,c+\frac{425}{1844}\,d \bigg)^2 + \frac{293135}{15072} \bigg( a\,d+\frac{215891}{293135}\,b\,d+\frac{215891}{293135}\,c \bigg)^2 + 12 \bigg( \frac{17}{24}\,a+b \bigg)^2 + \frac{10435033}{1172540} \bigg( b\,d+\frac{215891}{509026}\,c \bigg)^2 + \frac{29721597}{4072208} \bigg( c \bigg)^2 + \frac{287}{48} \bigg( a \bigg)^2 +  \frac{77709}{14752} \bigg( d \bigg)^2 + \frac{599}{408} \bigg( 1 \bigg)^2 \\
>0,
\end{dmath*}
\noindent \textit{Case 3}:
\begin{dmath*}
F_{2,1,1,1}(a,b,c)=42 + 42 \,a^2 + 180 \,a \,b + 42 \,b^2 + 102 \,a^2 \,b^2 + 180 \,a \,c + 180 \,b \,c + 420 \,a^2 \,b \,c + 420 \,a \,b^2 \,c + 42 \,c^2 + 102 \,a^2 \,c^2 + 420 \,a \,b \,c^2 + 102 \,b^2 \,c^2 + 942 \,a^2 \,b^2 \,c^2
\end{dmath*}
\begin{dmath*}
=942 \bigg( a\,b\,c+\frac{593}{3768}\,a+\frac{593}{3768}\,b+\frac{593}{3768}\,c \bigg)^2 + 102 \bigg( a\,b+\frac{247}{408}\,a\,c+\frac{247}{408}\,b\,c+\frac{109}{204} \bigg)^2 + \frac{105455}{1632} \bigg( a\,c+\frac{247}{655}\,b\,c+\frac{218}{655} \bigg)^2 +  \frac{72611}{1310} \bigg( b\,c+\frac{109}{451} \bigg)^2 + \frac{281375}{15072} \bigg( a+\frac{183407}{281375}\,b+\frac{183407}{281375}\,c \bigg)^2 + \frac{3021083}{281375} \bigg( b+\frac{183407}{464782}\,c \bigg)^2 + \frac{8426457}{929564} \bigg( c \bigg)^2 +  \frac{2241}{902} \bigg( 1 \bigg)^2\\
>0,
\end{dmath*}

\noindent where in each case, the first expression is obtained by (\ref{19B}) and the second expression (i.e. the SOS decomposition) is obtained by an application of \texttt{`SumsOfSquares'} to the first expression. \hfill \fbox

\begin{rem}
To establish the inequality (\ref{ineq1112}), we need only consider Case 1. However, to show that the equality sign holds if and only if $X_1,X_2,X_3,X_4$ are independent, we must check all three cases.
\end{rem}

\section{The SOS method with one exponent unbounded}\setcounter{equation}{0}

In this section, we demonstrate an even more powerful application of the SOS method to be used when one exponent is unknown. In particular, we extend Theorems \ref{thm234} and \ref{thm1112} to the case where $m_1$ is unbounded by exploiting the fact that verifying these GPIs can be reduced to proving the non-negativity of a multivariate polynomial with variables $a,b,...,m_1$.

\begin{thm}\label{thm34m}
Let $m \in \mathbb{N}$. For any centered Gaussian random vector $(X_1,X_2,X_3)$,
$$
E[X_1^{2m} X_2^{6} X_3^{4}]\ge E[X_1^{2m}] E[X_2^{6}] E[X_3^{4}].
$$
The equality holds if and only if $X_1,X_2,X_3$ are independent.
\end{thm}
\noindent{\bf Proof.}\ \
Let $X_3$ be a linear combination of $X_1$ and $X_2$, and let $U_1$ and $U_2$ be independent standard Gaussian random variables. Then, without loss of generality, we may write
$$
X_1=U_1, \,\, X_2=aU_1+U_2, \,\, X_3=bU_1+U_2, \,\,\, a,b \in \mathbb{R}.
$$
By Lemma \ref{lem2}, we need only show that $F_{m,3,1}$ and $F_{m,3,2}$ are strictly positive. By \cite[Theorem 3.5]{RusSun}, $F_{m,3,1}>0$. We will complete the proof by giving an SOS decomposition of $F_{m,3,2}$. By direct calculation (expanding and taking expectations), then letting $m=p^2+1$, we have that for any $p\in\mathbb{R}$,
\begin{dmath*}
\frac{F_{m,3,2}(a,b,p)}{2(2m-1)!!} = 450 + 2295 \,a^2 + 1620 \,a^4 + 135 \,a^6 + 3780 \,a \,b + 9000 \,a^3 \,b +
3780 \,a^5 \,b + 900 \,b^2 + 9990 \,a^2 \,b^2 + 14040 \,a^4 \,b^2 + 2790 \,a^6 \,b^2 +
 2700 \,a \,b^3 + 12600 \,a^3 \,b^3 + 11340 \,a^5 \,b^3 + 90 \,b^4 + 2295 \,a^2 \,b^4 +
 7020 \,a^4 \,b^4 + 5175 \,a^6 \,b^4 + 1575 \,a^2 \,p^2 + 1800 \,a^4 \,p^2 +
 213 \,a^6 \,p^2 + 2520 \,a \,b \,p^2 + 9600 \,a^3 \,b \,p^2 + 5112 \,a^5 \,b \,p^2 +
 630 \,b^2 \,p^2 + 10800 \,a^2 \,b^2 \,p^2 + 19170 \,a^4 \,b^2 \,p^2 +
 4464 \,a^6 \,b^2 \,p^2 + 2880 \,a \,b^3 \,p^2 + 17040 \,a^3 \,b^3 \,p^2 +
 17856 \,a^5 \,b^3 \,p^2 + 120 \,b^4 \,p^2 + \,3195 \,a^2 \,b^4 \,p^2 +
 11160 \,a^4 \,b^4 \,p^2 + 9129 \,a^6 \,b^4 \,p^2 + 450 \,a^4 \,p^4 + 90 \,a^6 \,p^4 +
 2400 \,a^3 \,b \,p^4 + 2160 \,a^5 \,b \,p^4 + 2700 \,a^2 \,b^2 \,p^4 +
 8100 \,a^4 \,b^2 \,p^4 + 2472 \,a^6 \,b^2 \,p^4 + 720 \,a \,b^3 \,p^4 +
 7200 \,a^3 \,b^3 \,p^4 + 9888 \,a^5 \,b^3 \,p^4 + 30 \,b^4 \,p^4 +
 1350 \,a^2 \,b^4 \,p^4 + 6180 \,a^4 \,b^4 \,p^4 + 6020 \,a^6 \,b^4 \,p^4 +
 12 \,a^6 \,p^6 + 288 \,a^5 \,b \,p^6 + 1080 \,a^4 \,b^2 \,p^6 + 576 \,a^6 \,b^2 \,p^6 +
 960 \,a^3 \,b^3 \,p^6 + 2304 \,a^5 \,b^3 \,p^6 + 180 \,a^2 \,b^4 \,p^6 +
 1440 \,a^4 \,b^4 \,p^6 + 1880 \,a^6 \,b^4 \,p^6 + 48 \,a^6 \,b^2 \,p^8 +
 192 \,a^5 \,b^3 \,p^8 + 120 \,a^4 \,b^4 \,p^8 + 280 \,a^6 \,b^4 \,p^8 + 16 \,a^6 \,b^4 \,p^{10}
\end{dmath*}
\begin{tiny}\begin{dmath*}
=\frac{2461933}{196} \bigg(-\frac{46403}{29543196}\,a^{3}b^{2}p^{5}+\frac{381808}{2461933}\,a^{3}b^{2}p^{3}+\frac{12180469}{29543196}\,a^{3}b^{2}p-\frac{8183}{4923866}\,a^{3}p^{3}+\frac{188501}{2461933}\,a^{2}b\,p^{3}+\frac{57820}{2461933}\,a\,b^{2}p^{3}+\frac{1233281}{29543196}\,a^{3}p+a^{2}b\,p+\frac{837655}{2461933}\,a\,b^{2}p+\frac{8357783}{39390928}\,a\,p+\frac{8216565}{91091521}\,b\,p \bigg)^2
+ 10899 \bigg(-\frac{2665}{87192}\,a^{3}b\,p^{3}+\frac{29641}{726600}\,a^{2}b^{2}p^{3}+\frac{233017}{523152}\,a^{3}b\,p+a^{2}b^{2}p+\frac{297779}{4272408}\,a^{2}p+\frac{1193}{4844}\,a\,b\,p-\frac{193}{7266}\,b^{2}p \bigg)^2
+ \frac{41695}{4} \bigg(\frac{107}{83390}\,a^{3}b^{2}p^{4}+\frac{37177}{500340}\,a^{3}b^{2}p^{2}+\frac{29011}{83390}\,a^{3}b^{2}-\frac{3259}{500340}\,a^{3}p^{2}+\frac{293859}{4086110}\,a^{2}b\,p^{2}-\frac{1574}{41695}\,a\,b^{2}p^{2}+\frac{118}{8339}\,a^{3}+a^{2}b+\frac{7844}{41695}\,a\,b^{2}+\frac{15677}{41695}\,a+\frac{2921}{41695}\,b \bigg)^2
+ \frac{70071801709219}{7090367040} \bigg(-\frac{245635624885}{70071801709219}\,a^{3}b^{2}p^{5}-\frac{11198510281680}{70071801709219}\,a^{3}b^{2}p^{3}+a^{3}b^{2}p+\frac{1830613217106}{350359008546095}\,a^{3}p^{3}-\frac{2274370721580}{70071801709219}\,a^{2}b\,p^{3}-\frac{1799591724552}{350359008546095}\,a\,b^{2}p^{3}-\frac{86820255707345}{1191220629056723}\,a^{3}p-\frac{13009331255460}{70071801709219}\,a\,b^{2}p-\frac{697867593895245}{13734073135006924}\,a\,p-\frac{81360974520540}{2592656663241103}\,b\,p \bigg)^2
+ \frac{172103}{20} \bigg(-\frac{420}{172103}\,a^{3}b\,p^{4}-\frac{3227}{1721030}\,a^{2}b^{2}p^{4}-\frac{161515}{2065236}\,a^{3}b\,p^{2}+\frac{16035}{172103}\,a^{2}b^{2}p^{2}+\frac{81745}{344206}\,a^{3}b+a^{2}b^{2}-\frac{208955}{16866094}\,a^{2}p^{2}+\frac{4780}{172103}\,a\,b\,p^{2}+\frac{1305}{344206}\,b^{2}p^{2}+\frac{10055}{172103}\,a^{2}+\frac{191845}{344206}\,a\,b-\frac{9953}{344206}\,b^{2}+\frac{10375}{172103} \bigg)^2
+ \frac{182859069551}{24016320} \bigg(-\frac{6418879954}{182859069551}\,a^{3}b^{2}p^{4}+a^{3}b^{2}p^{2}-\frac{41161024698}{182859069551}\,a^{3}b^{2}-\frac{20615021137}{914295347755}\,a^{3}p^{2}+\frac{154694373066}{1280013486857}\,a^{2}b\,p^{2}-\frac{40319988936}{914295347755}\,a\,b^{2}p^{2}+\frac{22432660800}{3108604182367}\,a^{3}-\frac{12775700256}{182859069551}\,a\,b^{2}+\frac{11681404428}{8960094407999}\,a-\frac{8345914044}{182859069551}\,b \bigg)^2
+ \frac{98097336461}{17210300} \bigg(-\frac{14512945525}{1177168037532}\,a^{3}b\,p^{4}+\frac{39929834765}{3335309439674}\,a^{2}b^{2}p^{4}+\frac{191874337775}{392389345844}\,a^{3}b\,p^{2}+a^{2}b^{2}p^{2}-\frac{221498760325}{1177168037532}\,a^{3}b+\frac{181601825675}{9613538973178}\,a^{2}p^{2}+\frac{24812335375}{196194672922}\,a\,b\,p^{2}-\frac{7224882875}{1667654719837}\,b^{2}p^{2}-\frac{111802291925}{4806769486589}\,a^{2}+\frac{57444211575}{196194672922}\,a\,b-\frac{403439350}{98097336461}\,b^{2}+\frac{267399874525}{3629601449057} \bigg)^2
+ \frac{1881745972159}{426892032} \bigg(-\frac{102378079626}{47043649303975}\,a^{3}b\,p^{3}-\frac{6247680732882}{47043649303975}\,a^{2}b^{2}p^{3}+a^{3}b\,p+\frac{40590877543374}{92205552635791}\,a^{2}p+\frac{371809874700}{13172221805113}\,a\,b\,p+\frac{35462399976}{1881745972159}\,b^{2}p \bigg)^2
+ \frac{1290111735374379}{365718139102} \bigg(-\frac{1254213486987157}{38703352061231370}\,a^{3}b^{2}p^{4}+a^{3}b^{2}+\frac{55648064894758831}{2193189950136444300}\,a^{3}p^{2}-\frac{813713883483332681}{7585857004001348520}\,a^{2}b\,p^{2}-\frac{919990297705031}{21501862256239650}\,a\,b^{2}p^{2}-\frac{16733411228329253}{219318995013644430}\,a^{3}-\frac{2332609290373289}{5160446941497516}\,a\,b^{2}-\frac{72218787593927536}{316077375166722855}\,a-\frac{5342404109573243}{25802234707487580}\,b \bigg)^2
+ \frac{281259138332925112160414689}{89209678367055858595200} \bigg(\frac{44836556278157635887192668}{281259138332925112160414689}\,a^{3}b^{2}p^{4}+\frac{2379918359966510057584675254}{23907026758298634533635248565}\,a^{3}p^{2}+a^{2}b\,p^{2}+\frac{512872285198992142407483276}{1406295691664625560802073445}\,a\,b^{2}p^{2}-\frac{246842735840239738068445884}{4781405351659726906727049713}\,a^{3}-\frac{41964775454104639080359730}{281259138332925112160414689}\,a\,b^{2}+\frac{82326511969366770611861466}{281259138332925112160414689}\,a+\frac{534358664600268384534896826}{10406588118318229149935343493}\,b \bigg)^2
+ \frac{724171982368916897}{282520329007680} \bigg(-\frac{78866109211794023}{3620859911844584485}\,a^{3}b\,p^{4}-\frac{169675800886557258}{12310923700271587249}\,a^{2}b^{2}p^{4}-\frac{169055740849158255}{12310923700271587249}\,a^{3}b\,p^{2}+a^{3}b-\frac{2958248375541023850}{35484427136076927953}\,a^{2}p^{2}-\frac{3538484271647279130}{35484427136076927953}\,a\,b\,p^{2}-\frac{36715685366745240}{12310923700271587249}\,b^{2}p^{2}+\frac{22124076308603404500}{35484427136076927953}\,a^{2}-\frac{49371151709178630}{724171982368916897}\,a\,b-\frac{1219853157005100}{724171982368916897}\,b^{2}-\frac{2201018396447410620}{26794363347649925189} \bigg)^2
+ \frac{26212746613691772978097}{11586751717902670352} \bigg(-\frac{117557911351656779877537}{6422122920354484379633765}\,a^{3}b\,p^{4}+\frac{109177880624035586572449}{2228083462163800703138245}\,a^{2}b^{2}p^{4}-\frac{28860678489162367020455027}{131011307575231481344528806}\,a^{3}b\,p^{2}-\frac{576983703856223757853069}{5137698336283587503707012}\,a^{2}p^{2}+\frac{7773801235190610828295310}{47523709610623184409289861}\,a\,b\,p^{2}-\frac{43531084034378371097507}{2228083462163800703138245}\,b^{2}p^{2}+\frac{51873196548430753987938}{183489226295842410846679}\,a^{2}+a\,b+\frac{2518705943463326354531}{26212746613691772978097}\,b^{2}+\frac{275045869188859060057562}{969871624706595600189589} \bigg)^2
+ \frac{53631965997981898754}{23881238132669869} \bigg(\frac{106353460685155026287}{429055727983855190032}\,a^{3}b\,p^{3}+\frac{4741131200225861740761}{10726393199596379750800}\,a^{2}b^{2}p^{3}+\frac{14207976389067404496925}{24027120767095890641792}\,a^{2}p+a\,b\,p+\frac{36998533915350406203}{268159829989909493770}\,b^{2}p \bigg)^2
+ \frac{1225720890883722229}{630646215382971}\bigg(a^{3}b^{2}p^{3}-\frac{227202054695158443}{12257208908837222290}\,a^{3}b^{2}p^{5}+\frac{1430298266691753873}{833490205800931115720}\,a^{3}p^{3}+\frac{249545321744126943}{19611534254139555664}\,a^{2}b\,p^{3}-\frac{14125355168335584543}{416745102900465557860}\,a\,b^{2}p^{3}-\frac{138247765820948296377}{1041862757251163894650}\,a^{3}p-\frac{106314906191237204067}{245144178176744445800}\,a\,b^{2}p-\frac{47876056145253870633}{960965178452838227536}\,a\,p-\frac{8296396344452300805}{90703345925395444946}\,b\,p \bigg)^2
+ \frac{9408995393797413663945824457325061}{5456620960508391197999624769900} \bigg(-\frac{350828535653402740933920328198325}{37635981575189654655783297829300244}\,a^{3}b\,p^{4}-\frac{2642069500129929412372761450964375}{37635981575189654655783297829300244}\,a^{2}b^{2}p^{4}+a^{3}b\,p^{2}+\frac{33326182700516389851833047354388525}{301087852601517237246266382634401952}\,a^{2}p^{2}-\frac{34382543630379550599089963715108825}{348132829570504305565995504921027257}\,a\,b\,p^{2}+\frac{309060196300930354417419284937655}{37635981575189654655783297829300244}\,b^{2}p^{2}+\frac{10118483817605990989094637303180925}{37635981575189654655783297829300244}\,a^{2}+\frac{23938506365002094892799773308875}{1017188691221342017723872914305412}\,b^{2}-\frac{454164830667551963943552068421383175}{5570125273128068889055928078736436112} \bigg)^2
+ \frac{8975651229529249369504496368487}{5625182766658502243208293780} \bigg(-\frac{7194504083747703440061891105659}{269269536885877481085134891054610}\,a^{3}b^{2}p^{4}+\frac{698906413627722036900956744291156}{37383587370989323623986227374748355}\,a^{3}p^{2}+\frac{22136513559812149638289245551787}{89756512295292493695044963684870}\,a\,b^{2}p^{2}-\frac{6691954163911113464596621358715}{152586070901997239281576438264279}\,a^{3}+a\,b^{2}+\frac{15404480886242227669327621841795}{3518455281975465752845762576446904}\,a+\frac{337263748895723350866231041032825}{664198190985164453343332731268038}\,b \bigg)^2
+ \frac{56993177287877132511208741}{49028835635348889160000}\bigg(a\,b^{2}p-\frac{86536462724224283284161220}{2906652041681733758071645791}\,a^{3}b^{2}p^{5}+\frac{91913835145634380748320645}{47475316680801651381836881253}\,a^{3}p^{3}-\frac{1942629801867344390693175}{16283764939393466431773926}\,a^{2}b\,p^{3}+\frac{94688574303636277492664890}{968884013893911252690548597}\,a\,b^{2}p^{3}-\frac{5849072067079347449891387284}{47475316680801651381836881253}\,a^{3}p-\frac{7244944468084291202534097925}{206657260845842482485642894866}\,a\,p+\frac{577737201220275968873509500}{2108747559651453902914723417}\,b\,p \bigg)^2
+ \frac{63640486372650027657933510345831052813}{59947138226161987319435290942729280} \bigg(\frac{6189698011362331106249631881668100060}{190921459117950082973800531037493158439}\,a^{3}b^{2}p^{5}+\frac{1613958574388828638236357715556086932}{63640486372650027657933510345831052813}\,a^{3}p^{3}+\frac{41873899779900945353977573965795898770}{63640486372650027657933510345831052813}\,a^{2}b\,p^{3}+\frac{318259569856788654338942826117002255662}{1081888268335050470184869675879127897821}\,a\,b^{2}p^{3}-\frac{2755819751181195952276013374464315408970}{14064547488355656112403305786428662671673}\,a^{3}p+a\,p+\frac{33729461638314340101874415820257108920}{63640486372650027657933510345831052813}\,b\,p \bigg)^2
+ \frac{219605904998627431295491364061}{255288158150393838069040000} \bigg(\frac{315464086652016426885241091137}{658817714995882293886474092183}\,a^{3}b\,p^{3}+a^{2}b^{2}p^{3}-\frac{2206498310122481015136534990475}{28695171586487317689277538237304}\,a^{2}p-\frac{46084543225252393769596876120}{219605904998627431295491364061}\,b^{2}p \bigg)^2
\end{dmath*}
\begin{dmath*}
+ \frac{751298069739776111324182015106815207080556113}{1051520590138866293781255650802756146623706} \bigg(\frac{1259576947571548888826894805716902024509294379}{10017307596530348150989093534757536094407414840}\,a^{3}b\,p^{4}+\frac{16691885447914782051714635337003715633373994653}{60103845579182088905934561208545216566444489040}\,a^{2}b^{2}p^{4}+\frac{6868670651718458722983859238852927535461291677}{16027692154448557041582549655612057751051863744}\,a^{2}p^{2}+a\,b\,p^{2}+\frac{538026084373165488785554207021731074975464139}{6010384557918208890593456120854521656644448904}\,b^{2}p^{2}-\frac{7046972944780067201238165098772217493260539121}{24041538231672835562373824483418086626577795616}\,a^{2}-\frac{2521771873821077005927768728351947837906418779}{30051922789591044452967280604272608283222244520}\,b^{2}+\frac{7946375997650306751578466826267626390558497827}{24041538231672835562373824483418086626577795616} \bigg)^2
+ \frac{132662645389010932413304395114633357161}{194974125872985827926077779809734912} \bigg(\frac{17215715386592930978709756023717434266434}{29849095212527459792993488900792505361225}\,a^{3}b\,p^{3}+a^{2}p-\frac{231023089838336117344966735522934217196}{1989939680835163986199565926719500357415}\,b^{2}p \bigg)^2
+ \frac{2372150908168012668643953810469261319417}{3940669915812521643187254085620532480}\bigg(a-\frac{722466858114816132429052833083361454628}{7116452724504038005931861431407783958251}\,a^{3}b^{2}p^{4}-\frac{178826903407899945980132343881925480318968}{3669717454935915598392196544795947261138099}\,a^{3}p^{2}+\frac{4010373078726940161739654321501343546996}{87769583602216468739826290987362668818429}\,a\,b^{2}p^{2}-\frac{7630499670216146862104053638396753924632}{40326565438856215366947214777977442430089}\,a^{3}+\frac{48125237977735494411552684017777899391380}{87769583602216468739826290987362668818429}\,b \bigg)^2
+ \frac{44632030280455886429070547400222520725609872436673}{96612369150139777965363789804339457701885721750}\bigg(a\,b^{2}p^{2}-\frac{76923411365460613986532612325215305933005020091713}{178528121121823545716282189600890082902439489746692}\,a^{3}b^{2}p^{4}+\frac{805039349079600695791544916353094327301134808218347}{28430603288650399655317938693941745702213488742160701}\,a^{3}p^{2}+\frac{3395782447168338760733203573364905836866229486755}{104141403987730401667831277267185881693089702352237}\,a^{3}+\frac{13339383660078040570890271971965271743159621388800}{44632030280455886429070547400222520725609872436673}\,b \bigg)^2
+ \frac{5798854105401307440575336297554933781480328659771088743}{18471594654158873019283056827099784316932251927685120} \bigg(\frac{30277127948142456088437246825356204288043856074207028}{5798854105401307440575336297554933781480328659771088743}\,a^{3}b\,p^{4}+\frac{5782049733604745109787099157079690598436073253290462402}{98580519791822226489780717058433874285165587216108508631}\,a^{2}b^{2}p^{4}-\frac{5250630188784926838577066099329704582123901070477120885}{150770206740433993454958743736428278318488545154048307318}\,a^{2}p^{2}+\frac{168217991306225563253703841884035139511111244349194852956}{3647479232297422380121886531162053348551126726996014819347}\,b^{2}p^{2}+a^{2}-\frac{286373275422327888484000557331418004478397041314559612}{5798854105401307440575336297554933781480328659771088743}\,b^{2}-\frac{55402541201255288223393347052711069475278386101246562353}{214557601899848375301287443009532549914772160411530283491} \bigg)^2
+ \frac{6023346958102712039190062093270040252063816812126410557433925809861}{28913605404605437364487949936098114669496215305584495511450013696}\bigg(a^{2}p^{2}+\frac{1932803082016863792998416665346552216422730977358524499321554587736}{6023346958102712039190062093270040252063816812126410557433925809861}\,a^{3}b\,p^{4}-\frac{1966908786241381431097115740344442503599505523628009324880191741608}{30116734790513560195950310466350201260319084060632052787169629049305}\,a^{2}b^{2}p^{4}+\frac{1631333850223511517336093223712347313658884881839050334805943101656}{187558675081515142210423220627072045472680236377599615872571254178345}\,b^{2}p^{2}+\frac{164250264269349073981479303374507844208541859793506504207049023090168}{1114319187249001727250161487254957446631806110243385953125276274824285}\,b^{2}-\frac{414016067374660616574259918319937284482476074583081467069433565229398}{1114319187249001727250161487254957446631806110243385953125276274824285} \bigg)^2
+ \frac{20417388528673696915508734078234521590979868225263}{119407108045100265632312872315209014648092782928} \bigg(\frac{2787395580558409915196547884788149626523493488602}{183756496758063272239578606704110694318818814027367}\,a^{3}b^{2}p^{5}-\frac{772960345519102603146661898113933611086632403522}{6805796176224565638502911359411507196993289408421}\,a^{3}p^{3}-\frac{1950870272617726351443232491042466090559368814603}{20417388528673696915508734078234521590979868225263}\,a^{2}b\,p^{3}+\frac{1153578369605926458091217357612762660459762317997}{34028980881122828192514556797057535984966447042105}\,a\,b^{2}p^{3}+a^{3}p+\frac{2019305919977093630271343170846884619431723059891579}{15108867511218535717476463217893545977325102486694620}\,b\,p \bigg)^2
+ \frac{97584090849321545201161697823073171144143547991904113621}{606995611814200055435359444643026281868294265138752800} \bigg(a^{3}b^{2}p^{4}-\frac{31810482542365387581229673409105191394901688694296195918080}{62161065871017824293140001513297610018819440070842920376577}\,a^{3}p^{2}+\frac{100601823575352792714848746666196094406703716096722184464}{683088635945250816408131884761512198009004835943328795347}\,a^{3}+\frac{587538523063443535832416944491098895488438015055366215340}{3610611361424897172442982819453707332333311275700452203977}\,b \bigg)^2
+ \frac{2547979952487198785212709307211844770932994888524125459}{16007232606480178381758847517335864927328216688606192} \bigg(\frac{62924382969983534282796669999402329145434355225981759476}{389840932730541414137544524003412249952748217944191195227}\,a^{3}b^{2}p^{5}+\frac{22014889552078903961755702990501698619139381602045081767}{216578295961411896743080291113006805529304565524550664015}\,a^{3}p^{3}+a^{2}b\,p^{3}+\frac{87275259705742649259330286073835161560414133380363932633}{433156591922823793486160582226013611058609131049101328030}\,a\,b^{2}p^{3}-\frac{9403912318199261971375257692393709136891230429950032823483}{32053587802288960717975883084725007218337075697633498274220}\,b\,p \bigg)^2
+ \frac{1099497882555389016616929857482383397985239609783}{7307058508026722157324806082914005312427880000} \bigg(a^{3}b\,p^{3}+\frac{8590683514112276150746077669097583411718303463050}{40681421654549393614826404726848185725453865561971}\,b^{2}p \bigg)^2
+ \frac{87093407027768187977750098712984102927291902454731023702879071183}{806468269105590251664273218411681181613360824552458816579375200}\bigg(b\,p-\frac{83125712946748370121615793005973132886828547335592924208551617840}{783840663249913691799750888416856926345627122092579213325911640647}\,a^{3}b^{2}p^{5}+\frac{7080913565140181275242944969423540790010775212419027490041319481}{174186814055536375955500197425968205854583804909462047405758142366}\,a^{3}p^{3}+\frac{270622143972119826134314027023228971339653128390218279989585760815}{348373628111072751911000394851936411709167609818924094811516284732}\,a\,b^{2}p^{3} \bigg)^2
+ \frac{36264290045723190857537954223173726795432470564619545082170843587}{483660908681034840557777781005388711800275874089673953545443350}\bigg(a^{3}-\frac{7875000105462815664557044530875488162146430465438705584698594227565}{19800302364964862208215723005852854830306128928282271614865280598502}\,a^{3}p^{2}+\frac{1026777627345518238625370370074903889261048605992655266995572575265}{8050672390150548370373425837544567348586008465345539008241927276314}\,b \bigg)^2
\end{dmath*}
\begin{dmath*}
+ \frac{46124565169902990116262869292850250792875255601422953862150214947411658251373}{824352720463386196367795989248493895120777369755011430114035281352837256680} \bigg(a^{3}p^{2}+\frac{90943024873640246533155452571824687490347492900547725559213577307263460482883}{368996521359223920930102954342802006343002044811383630897201719579293266010984}\,b \bigg)^2 + \frac{11654614620727252139843764945642431404074405212468642468332588577927793423}{217593408861460472415740993119380204105805382338066581387300569881228625}\bigg(a^{2}b^{2}p^{4}+\frac{32479559436844990478561515462508498134705303006084848125725032744406011155}{69927687724363512839062589673854588424446431274811854809995531467566760538}\,a^{3}b\,p^{4}-\frac{406683476176259121096934802645185572809960567684558370081042577411740860625}{862441481933816658348438605977539923901505985722679542656611554766656713302}\,b^{2}p^{2}+\frac{949216098171298274411142425092007893258330674494823501593205082774276676585}{3449765927735266633393754423910159695606023942890718170626446219066626853208}\,b^{2}-\frac{2932964600106421094570667692811376749159303863835736746419697511042555299635}{13799063710941066533575017695640638782424095771562872682505784876266507412832} \bigg)^2
+ \frac{986460215342179026424390748255258760026853906327490295562098667618081992989242189325219}{18913049541014711387883961761903880243832853287543600099743197129655724894968748874240} \bigg(b \bigg)^2
+ \frac{1218803678994326257684500778213100903032931005357009011}{25588614220711568583725808573187508821310481438479759} \bigg(b^{2}p \bigg)^2
+ \frac{19166264940936963008123131317049546099244670234548632987602016878308827520858144607}{408452285843855569393820523790962907959753234838261031402171232337488619419827200}\bigg(1+\frac{5843119416744777965202677384743069748373864488455792264247623166257485440193041200}{57498794822810889024369393951148638297734010703645898962806050634926482562574433821}\,a^{3}b\,p^{4}-\frac{92925645164886570099270193922892054314034811758914630736667246371180720076870377424}{325826503995928371138093232389842283687159393987326760789234286931250067854588458319}\,b^{2}p^{2}-\frac{9082038525278957153218384242323589269424217122856914454603945973297625726562262612}{19166264940936963008123131317049546099244670234548632987602016878308827520858144607}\,b^{2} \bigg)^2
+ \frac{23578567916368233618521083347021261517635043797485493835603740546535154013719788691954786}{655965417603567558953014169326020715246648838777426964000679027660119621901369999174575}\bigg(b^{2}-\frac{34667926262155625075188570818536417330265957461567036799016066419955816630658705711640977}{565885629992837606844506000328510276423241051139651852054489773116843696329274928606914864}\,a^{3}b\,p^{4}-\frac{37979020513911577022337194238075360554905610397228550555552444549892004248442125806893445}{246668095125083059393759025784222428184489688958309781664777593409906226605068558623526992}\,b^{2}p^{2} \bigg)^2
+ \frac{25985558106761789167791401489470178468552050845136456244416868845005311497}{1610879656385600404836465825795353967743191027802705014408451300600768000}\bigg(a\,b^{2}p^{3}-\frac{101430110962951336045493287382728911174877170403422999192646604922231035200}{233870022960856102510122613405231606216968457606228106199751819605047803473}\,a^{3}b^{2}p^{5}+\frac{76297795422328098743873832844677774083208815967704333819594992030941205566}{1273292347231327669221778672984038744959050491411686355976426573405260263353}\,a^{3}p^{3} \bigg)^2
+ \frac{1247026159056318920889118972193317004988891724553734109809904660526782706729276279}{180310929291428311238496077881269726673651140088808904869821767059918905893418330}\bigg(a^{3}p^{3}-\frac{13073485100002825609349930742072279113522220600160995746843509216184631212540008159}{22446470863013740576004141499479706089800051041967213976578283889482088721126973022}\,a^{3}b^{2}p^{5} \bigg)^2
+ \frac{6635364387713420500497167118347376568266566032959160711499775776758948454187251786344377759541093}{1148141315569211608254190761415241714227525706225348709736873786285749522355952112969068736963200}\bigg(b^{2}p^{2}+\frac{11033434426734849202372000295956059607726548933570936180092623498588194212727967314509876582015657}{19906093163140261501491501355042129704799698098877482134499327330276845362561755359033133278623279}\,a^{3}b\,p^{4} \bigg)^2
+ \frac{388681935647012442251946102542793082043360840320090931730502746670863314264566630906737289084953671697}{172519474080548933012926345077031790774930716856938178498994170195732659808868546444953821748068418000} \bigg(a^{3}b\,p^{4} \bigg)^2
+ \frac{2629929480562269930314071400782243373877386514044892496969174184208223685867430351966437}{1167665414293974784763735440802934310791398655203134471061602327930858255273025136604440} \bigg(a^{3}b^{2}p^{5} \bigg)^2\\
>0,
\end{dmath*}\end{tiny}

\noindent where the second expression (i.e. the SOS decomposition) is obtained by an application of \texttt{`SumsOfSquares'} to the first expression. \hfill \fbox

\begin{thm}\label{thm111m}
Let $m\in\mathbb{N}$. For any centered Gaussian random vector $(X_1,X_2,X_3,X_4)$,
\begin{equation}\label{ineq111m}
E[X_1^{2m} X_2^{2} X_3^{2} X_4^{2}]\ge E[X_1^{2m}] E[X_2^{2}] E[X_3^{2}] E[X_4^{2}].
\end{equation}
The equality holds if and only if $X_1,X_2,X_3,X_4$ are independent.
\end{thm}
\noindent{\bf Proof.}\ \
Let $X_4$ be a linear combination of $X_1$, $X_2$ and $X_3$, and let $U_1$, $U_2$ and $U_3$ be independent standard Gaussian random variables. Then, without loss of generality, we may consider 3 constructions of $X_1,\dots,X_4$ in terms of $U_1,U_2,U_3$:

\noindent \textit{Case 1}:
$$
X_1=U_1, \,\, X_2=aU_1+U_2, \,\, X_3=bU_1+cU_2+U_3, \,\, X_4=dU_1+eU_2+U_3, \,\,\, a,b,c,d,e \in \mathbb{R}.
$$
\textit{Case 2}:
$$
X_1=U_1, \,\, X_2=aU_1+U_2, \,\, X_3=bU_1+U_2, \,\, X_4=cU_1+dU_2+U_3, \,\,\, a,b,c,d \in \mathbb{R}.
$$
\textit{Case 3}:
$$
X_1=U_1, \,\, X_2=aU_1+U_2, \,\, X_3=bU_1+U_2, \,\, X_4=cU_1+U_2, \,\,\, a,b,c \in \mathbb{R}.
$$
By Lemma \ref{lem2}, we need only show that, in each of these $3$ cases, $F_{m,1,1,1}>0$. By direct calculation (expanding and taking expectations), then letting $m=p^2+1$, we have that for any $p\in\mathbb{R}$,

\noindent \textit{Case 1}:
\begin{dmath*}
\frac{F_{m,1,1,1}(a,b,c,d,e,p)}{2(2m-1)!!}=1 + 4 \,a^2 + b^2 + 7 \,a^2 \,b^2 + 6 \,a \,b \,c + c^2 + a^2 \,c^2 + 6 \,b \,d +
 30 \,a^2 \,b \,d + 12 \,a \,c \,d + d^2 + 7 \,a^2 \,d^2 + 7 \,b^2 \,d^2 +
 52 \,a^2 \,b^2 \,d^2 + 30 \,a \,b \,c \,d^2 + 4 \,c^2 \,d^2 + 7 \,a^2 \,c^2 \,d^2 +
 12 \,a \,b \,e + 6 \,c \,e + 6 \,a^2 \,c \,e + 6 \,a \,d \,e + 30 \,a \,b^2 \,d \,e + 18 \,b \,c \,d \,e +
 30 \,a^2 \,b \,c \,d \,e + 18 \,a \,c^2 \,d \,e + e^2 + a^2 \,e^2 + 4 \,b^2 \,e^2 +
 7 \,a^2 \,b^2 \,e^2 + 18 \,a \,b \,c \,e^2 + 7 \,c^2 \,e^2 + 4 \,a^2 \,c^2 \,e^2 +
 3 \,a^2 \,p^2 + b^2 \,p^2 + 8 \,a^2 \,b^2 \,p^2 + 4 \,a \,b \,c \,p^2 + a^2 \,c^2 \,p^2 +
 4 \,b \,d \,p^2 + 32 \,a^2 \,b \,d \,p^2 + 8 \,a \,c \,d \,p^2 + d^2 \,p^2 + 8 \,a^2 \,d^2 \,p^2 +
 8 \,b^2 \,d^2 \,p^2 + 71 \,a^2 \,b^2 \,d^2 \,p^2 + 32 \,a \,b \,c \,d^2 \,p^2 +
 3 \,c^2 \,d^2 \,p^2 + 8 \,a^2 \,c^2 \,d^2 \,p^2 + 8 \,a \,b \,e \,p^2 + 4 \,a^2 \,c \,e \,p^2 +
 4 \,a \,d \,e \,p^2 + 32 \,a \,b^2 \,d \,e \,p^2 + 12 \,b \,c \,d \,e \,p^2 +
 32 \,a^2 \,b \,c \,d \,e \,p^2 + 12 \,a \,c^2 \,d \,e \,p^2 + a^2 \,e^2 \,p^2 +
 3 \,b^2 \,e^2 \,p^2 + 8 \,a^2 \,b^2 \,e^2 \,p^2 + 12 \,a \,b \,c \,e^2 \,p^2 +
 3 \,a^2 \,c^2 \,e^2 \,p^2 + 2 \,a^2 \,b^2 \,p^4 + 8 \,a^2 \,b \,d \,p^4 + 2 \,a^2 \,d^2 \,p^4 +
 2 \,b^2 \,d^2 \,p^4 + 30 \,a^2 \,b^2 \,d^2 \,p^4 + 8 \,a \,b \,c \,d^2 \,p^4 +
 2 \,a^2 \,c^2 \,d^2 \,p^4 + 8 \,a \,b^2 \,d \,e \,p^4 + 8 \,a^2 \,b \,c \,d \,e \,p^4 +
 2 \,a^2 \,b^2 \,e^2 \,p^4 + 4 \,a^2 \,b^2 \,d^2 \,p^6
\end{dmath*}
\begin{tiny}\begin{dmath*}
=\frac{123}{2}\bigg(\frac{27}{410}\,a\,b\,d\,p^{3}+a\,b\,d\,p+\frac{137}{1230}\,a\,c\,e\,p+\frac{71}{615}\,c\,d\,p+\frac{71}{615}\,b\,e\,p+\frac{9}{82}\,a\,p \bigg)^2
+52\bigg(\frac{19}{208}\,a\,b\,d\,p^{2}+a\,b\,d+\frac{83}{416}\,a\,c\,e+\frac{81}{416}\,c\,d+\frac{81}{416}\,b\,e+\frac{85}{416}\,a \bigg)^2
+\frac{89299}{4160}\bigg(a\,b\,d\,p^{2}+\frac{3301}{25514}\,a\,c\,e+\frac{20177}{178598}\,c\,d+\frac{20177}{178598}\,b\,e+\frac{3155}{25514}\,a \bigg)^2
+7\bigg(\frac{5}{14}\,a\,b\,p^{2}+\frac{5}{14}\,a\,d\,p^{2}+a\,b+\frac{5}{8}\,a\,d+\frac{9}{28}\,c+\frac{55}{224}\,e \bigg)^2
+7\bigg(\frac{13}{35}\,a\,c\,d\,p^{2}+\frac{9}{28}\,a\,b\,e\,p^{2}+\frac{13}{35}\,b\,d\,p^{2}+\frac{37}{56}\,a\,c\,d+a\,b\,e+\frac{39}{56}\,b\,d+\frac{45}{56}\,c\,e+\frac{59}{224} \bigg)^2
+\frac{273}{64}\bigg(\frac{20}{91}\,a\,b\,p^{2}+\frac{20}{91}\,a\,d\,p^{2}+a\,d+\frac{20}{273}\,c+\frac{43}{156}\,e \bigg)^2
+\frac{1767}{448}\bigg(\frac{1192}{8835}\,a\,c\,d\,p^{2}+\frac{2494}{8835}\,a\,b\,e\,p^{2}+\frac{104}{465}\,b\,d\,p^{2}+a\,c\,d+\frac{13}{31}\,b\,d+\frac{15}{31}\,c\,e+\frac{59}{372} \bigg)^2
+\frac{30613}{8200}\bigg(a\,b\,d\,p^{3}+\frac{15776}{30613}\,a\,c\,e\,p+\frac{15641}{30613}\,c\,d\,p+\frac{15641}{30613}\,b\,e\,p+\frac{16855}{30613}\,a\,p \bigg)^2
+\frac{7}{2}\bigg(\frac{9}{140}\,a\,c\,d\,p+a\,b\,e\,p+\frac{1}{10}\,b\,d\,p \bigg)^2
+\frac{19519}{5600}\bigg(a\,c\,d\,p+\frac{14}{149}\,b\,d\,p \bigg)^2
+3\bigg(a\,b\,p+\frac{5}{24}\,a\,d\,p \bigg)^2
+\frac{361}{124}\bigg(\frac{351}{1805}\,a\,c\,d\,p^{2}+\frac{351}{1805}\,a\,b\,e\,p^{2}+\frac{353}{1805}\,b\,d\,p^{2}+b\,d+\frac{287}{722}\,c\,e+\frac{411}{1444} \bigg)^2
+\frac{551}{192}\bigg(a\,d\,p \bigg)^2
+\frac{3627}{1490}\bigg(b\,d\,p \bigg)^2
+\frac{6267251}{3571960}\bigg(\frac{7437577}{12534502}\,a\,c\,e+\frac{3141786}{6267251}\,c\,d+b\,e-\frac{3928365}{25069004}\,a \bigg)^2
+\frac{1169378377}{802208128}\bigg(\frac{260572758}{1169378377}\,a\,c\,e-\frac{109994220}{1169378377}\,c\,d+a \bigg)^2
+\frac{6084238559}{4677513508}\bigg(\frac{367352216}{869176937}\,a\,c\,e+c\,d \bigg)^2
+\frac{915139}{734712}\bigg(a\,c\,e\,p+\frac{717236}{915139}\,c\,d\,p+\frac{717236}{915139}\,b\,e\,p-\frac{136728}{915139}\,a\,p \bigg)^2
+\frac{1006310}{915139}\bigg(a\,p-\frac{8122597}{32201920}\,c\,d\,p-\frac{8122597}{32201920}\,b\,e\,p \bigg)^2
+\frac{6333}{5776}\bigg(\frac{19586}{31665}\,a\,c\,d\,p^{2}+\frac{19586}{31665}\,a\,b\,e\,p^{2}+\frac{4886}{10555}\,b\,d\,p^{2}+c\,e+\frac{3019}{12666} \bigg)^2
+1\bigg(a\,c+\frac{3}{8}\,a\,e+\frac{3}{4}\,b+\frac{11}{32}\,d \bigg)^2
+1\bigg(\frac{3}{8}\,a\,c\,p+a\,e\,p+\frac{17}{32}\,b\,p+\frac{7}{8}\,d\,p \bigg)^2
+\frac{82}{91}\bigg(\frac{73}{164}\,a\,b\,p^{2}+a\,d\,p^{2}+\frac{463}{2624}\,c+\frac{23}{82}\,e \bigg)^2
+\frac{903}{1024}\bigg(\frac{212}{301}\,a\,e+\frac{376}{903}\,b+d \bigg)^2
+\frac{55}{64}\bigg(a\,c\,p+\frac{173}{220}\,b\,p+\frac{13}{55}\,d\,p \bigg)^2
+\frac{237}{328}\bigg(a\,b\,p^{2}+\frac{955}{3792}\,c+\frac{61}{948}\,e \bigg)^2
+\frac{2253517443}{3476707748}\bigg(a\,c\,e \bigg)^2
+\frac{77689}{158325}\bigg(\frac{146437}{621512}\,a\,c\,d\,p^{2}+\frac{146437}{621512}\,a\,b\,e\,p^{2}+b\,d\,p^{2}+\frac{22520}{77689} \bigg)^2
+\frac{127}{301}\bigg(a\,e-\frac{1889}{4064}\,b \bigg)^2
+\frac{114972691}{283409472}\bigg(\frac{582802531}{2874317275}\,a\,c\,d\,p^{2}+a\,b\,e\,p^{2}-\frac{133581818}{574863455} \bigg)^2
+\frac{111810532673}{287431727500}\bigg(a\,c\,d\,p^{2}-\frac{17576555}{90976837} \bigg)^2
+\frac{381989069}{1030461440}\bigg(c\,d\,p+\frac{124373709}{381989069}\,b\,e\,p \bigg)^2
+\frac{253181389}{763978138}\bigg(b\,e\,p \bigg)^2
+\frac{75473}{390144}\bigg(b \bigg)^2
+\frac{41}{220}\bigg(b\,p+\frac{883}{1312}\,d\,p \bigg)^2
+\frac{38225}{212352}\bigg(c-\frac{59373}{152900}\,e \bigg)^2
+\frac{35901007}{234854400}\bigg(e \bigg)^2
+\frac{17121}{167936}\bigg(d\,p \bigg)^2
+\frac{229159727}{5822517568}\bigg(1 \bigg)^2\\
>0,
\end{dmath*}\end{tiny}

\noindent \textit{Case 2}:
\begin{dmath*}
\frac{F_{m,1,1,1}(a,b,c,d,p)}{2(2m-1)!!}=1+a^2+6\,a\,b+b^2+7\,a^2\,b^2+4\,c^2+7\,a^2\,c^2+30\,a\,b\,c^2+7\,b^2\,c^2+52\,a^2\,b^2\,c^2+18\,a\,c\,d+18\,b\,c\,d+30\,a^2\,b\,c\,d+30\,a\,b^2\,c\,d+7\,d^2+4\,a^2\,d^2+18\,a\,b\,d^2+4\,b^2\,d^2+7\,a^2\,b^2\,d^2+a^2\,p^2+4\,a\,b\,p^2+b^2\,p^2+8\,a^2\,b^2\,p^2+3\,c^2\,p^2+8\,a^2\,c^2\,p^2+32\,a\,b\,c^2\,p^2+8\,b^2\,c^2\,p^2+71\,a^2\,b^2\,c^2\,p^2+12\,a\,c\,d\,p^2+12\,b\,c\,d\,p^2+32\,a^2\,b\,c\,d\,p^2+32\,a\,b^2\,c\,d\,p^2+3\,a^2\,d^2\,p^2+12\,a\,b\,d^2\,p^2+3\,b^2\,d^2\,p^2+8\,a^2\,b^2\,d^2\,p^2+2\,a^2\,b^2\,p^4+2\,a^2\,c^2\,p^4+8\,a\,b\,c^2\,p^4+2\,b^2\,c^2\,p^4+30\,a^2\,b^2\,c^2\,p^4+8\,a^2\,b\,c\,d\,p^4+8\,a\,b^2\,c\,d\,p^4+2\,a^2\,b^2\,d^2\,p^4+4\,a^2\,b^2\,c^2\,p^6
\end{dmath*}
\begin{small}\begin{dmath*}
=\frac{368}{5}\bigg(\frac{25}{2944}\,a\,b\,c\,p^{3}+a\,b\,c\,p+\frac{281}{2944}\,a\,d\,p+\frac{281}{2944}\,b\,d\,p+\frac{281}{2944}\,c\,p \bigg)^2
+52\bigg(\frac{41}{208}\,c-\frac{1}{40}\,a\,b\,c\,p^{2}+a\,b\,c+\frac{41}{208}\,a\,d+\frac{41}{208}\,b\,d \bigg)^2
+\frac{11487}{400}\bigg(a\,b\,c\,p^{2}+\frac{2825}{22974}\,a\,d+\frac{2825}{22974}\,b\,d+\frac{2825}{22974}\,c \bigg)^2
+7\bigg(\frac{73}{280}\,a\,b\,p^{2}+a\,b+\frac{9}{28} \bigg)^2
+7\bigg(\frac{19}{56}\,a\,b\,d\,p^{2}+\frac{53}{140}\,a\,c\,p^{2}+\frac{53}{140}\,b\,c\,p^{2}+a\,b\,d+\frac{19}{28}\,a\,c+\frac{19}{28}\,b\,c+\frac{23}{28}\,d \bigg)^2
+\frac{87}{20}\bigg(a\,b\,p \bigg)^2
+\frac{94083}{23552}\bigg(a\,b\,c\,p^{3}+\frac{18177}{31361}\,a\,d\,p+\frac{18177}{31361}\,b\,d\,p+\frac{18177}{31361}\,c\,p \bigg)^2
+\frac{423}{112}\bigg(\frac{1163}{4230}\,a\,b\,d\,p^{2}+\frac{323}{2115}\,a\,c\,p^{2}+\frac{53}{235}\,b\,c\,p^{2}+a\,c+\frac{19}{47}\,b\,c+\frac{23}{47}\,d \bigg)^2
+\frac{13}{4}\bigg(a\,b\,d\,p+\frac{8}{65}\,a\,c\,p+\frac{8}{65}\,b\,c\,p \bigg)^2
+\frac{4161}{1300}\bigg(a\,c\,p+\frac{8}{73}\,b\,c\,p \bigg)^2
+\frac{4617}{1460}\bigg(b\,c\,p \bigg)^2
+\frac{297}{94}\bigg(\frac{1163}{5940}\,a\,b\,d\,p^{2}+\frac{1163}{5940}\,a\,c\,p^{2}+\frac{437}{5940}\,b\,c\,p^{2}+b\,c+\frac{23}{66}\,d \bigg)^2
+\frac{3692291}{2389296}\bigg(a\,d+\frac{1900319}{3692291}\,b\,d+\frac{1900319}{3692291}\,c \bigg)^2
+\frac{17071}{11200}\bigg(a\,b\,p^{2}+\frac{6030}{17071} \bigg)^2
+\frac{8388915}{7384582}\bigg(b\,d+\frac{1900319}{5592610}\,c \bigg)^2
+\frac{22478787}{22370440}\bigg(c \bigg)^2
+1\bigg(a\,p+\frac{7}{8}\,b\,p \bigg)^2
+1\bigg(\frac{3}{4}\,a+b \bigg)^2
+\frac{87}{88}\bigg(\frac{1519}{2610}\,a\,b\,d\,p^{2}+\frac{1519}{2610}\,a\,c\,p^{2}+\frac{1519}{2610}\,b\,c\,p^{2}+d \bigg)^2
+\frac{619373}{627220}\bigg(\frac{34867}{141790}\,a\,b\,d\,p^{2}+\frac{34867}{141790}\,a\,c\,p^{2}+b\,c\,p^{2} \bigg)^2
+\frac{14179}{31320}\bigg(a \bigg)^2
+\frac{7}{16}\bigg(\frac{462568}{619373}\,a\,d\,p+b\,d\,p+\frac{462568}{619373}\,c\,p \bigg)^2
+\frac{1081941}{2477492}\bigg(\frac{462568}{1081941}\,a\,d\,p+c\,p \bigg)^2
+\frac{217111453}{510444000}\bigg(a\,b\,d\,p^{2}+\frac{34867}{176657}\,a\,c\,p^{2} \bigg)^2
+\frac{21663583}{52997100}\bigg(a\,c\,p^{2} \bigg)^2
+\frac{1544509}{4327764}\bigg(a\,d\,p \bigg)^2
+\frac{15}{64}\bigg(b\,p \bigg)^2
+\frac{2957}{34142}\bigg(1 \bigg)^2\\
>0,
\end{dmath*}\end{small}
\noindent \textit{Case 3}:
\begin{dmath*}
\frac{F_{m,1,1,1}(a,b,c,p)}{2(2m-1)!!}=7+4\,a^2+18\,a\,b+4\,b^2+7\,a^2\,b^2+18\,a\,c+18\,b\,c+30\,a^2\,b\,c+30\,a\,b^2\,c+4\,c^2+7\,a^2\,c^2+30\,a\,b\,c^2+7\,b^2\,c^2+52\,a^2\,b^2\,c^2+3\,a^2\,p^2+12\,a\,b\,p^2+3\,b^2\,p^2+8\,a^2\,b^2\,p^2+12\,a\,c\,p^2+12\,b\,c\,p^2+32\,a^2\,b\,c\,p^2+32\,a\,b^2\,c\,p^2+3\,c^2\,p^2+8\,a^2\,c^2\,p^2+32\,a\,b\,c^2\,p^2+8\,b^2\,c^2\,p^2+71\,a^2\,b^2\,c^2\,p^2+2\,a^2\,b^2\,p^4+8\,a^2\,b\,c\,p^4+8\,a\,b^2\,c\,p^4+2\,a^2\,c^2\,p^4+8\,a\,b\,c^2\,p^4+2\,b^2\,c^2\,p^4+30\,a^2\,b^2\,c^2\,p^4+4\,a^2\,b^2\,c^2\,p^6
\end{dmath*}
\begin{small}
\begin{dmath*}
=\frac{372}{5}\bigg(\frac{5}{2976}\,a\,b\,c\,p^{3}+a\,b\,c\,p+\frac{139}{1488}\,a\,p+\frac{139}{1488}\,b\,p+\frac{139}{1488}\,c\,p \bigg)^2
+52\bigg(\frac{81}{416}\,c-\frac{17}{520}\,a\,b\,c\,p^{2}+a\,b\,c+\frac{81}{416}\,a+\frac{81}{416}\,b \bigg)^2
+\frac{154411}{5200}\bigg(a\,b\,c\,p^{2}+\frac{76045}{617644}\,a+\frac{76045}{617644}\,b+\frac{76045}{617644}\,c \bigg)^2
+7\bigg(\frac{13}{28}\,a\,b\,p^{2}+\frac{13}{28}\,a\,c\,p^{2}+\frac{13}{28}\,b\,c\,p^{2}+\frac{23}{28}\,a\,b+\frac{23}{28}\,a\,c+\frac{23}{28}\,b\,c+1 \bigg)^2
+\frac{95227}{23808}\bigg(a\,b\,c\,p^{3}+\frac{56266}{95227}\,a\,p+\frac{56266}{95227}\,b\,p+\frac{56266}{95227}\,c\,p \bigg)^2
+\frac{63}{20}\bigg(a\,b\,p+\frac{13}{126}\,a\,c\,p+\frac{13}{126}\,b\,c\,p \bigg)^2
+\frac{15707}{5040}\bigg(a\,c\,p+\frac{13}{139}\,b\,c\,p \bigg)^2
+\frac{2147}{695}\bigg(b\,c\,p \bigg)^2
+\frac{255}{112}\bigg(\frac{1}{75}\,a\,b\,p^{2}-\frac{137}{1275}\,a\,c\,p^{2}+\frac{1}{75}\,b\,c\,p^{2}+\frac{1}{15}\,a\,b+a\,c+\frac{1}{15}\,b\,c \bigg)^2
+\frac{34}{15}\bigg(\frac{1}{80}\,a\,b\,p^{2}+\frac{7}{340}\,a\,c\,p^{2}-\frac{37}{340}\,b\,c\,p^{2}+\frac{1}{16}\,a\,b+b\,c \bigg)^2
+\frac{289}{128}\bigg(a\,b-\frac{159}{1445}\,a\,b\,p^{2}+\frac{28}{1445}\,a\,c\,p^{2}+\frac{28}{1445}\,b\,c\,p^{2} \bigg)^2
+\frac{15598365}{9882304}\bigg(\frac{2728879}{5199455}\,a+b+\frac{2728879}{5199455}\,c \bigg)^2
+\frac{11892501}{10398910}\bigg(a+\frac{2728879}{7928334}\,c \bigg)^2
+\frac{10657213}{10571112}\bigg(c \bigg)^2
+\frac{908823}{952270}\bigg(a\,p+\frac{1341511}{1817646}\,b\,p+\frac{1341511}{1817646}\,c\,p \bigg)^2
+\frac{669}{1445}\bigg(a\,b\,p^{2}+\frac{3571}{13380}\,a\,c\,p^{2}+\frac{3571}{13380}\,b\,c\,p^{2} \bigg)^2
+\frac{3159157}{7270584}\bigg(b\,p+\frac{1341511}{3159157}\,c\,p \bigg)^2
+\frac{9780727}{22746000}\bigg(\frac{3571}{16951}\,a\,c\,p^{2}+b\,c\,p^{2} \bigg)^2
+\frac{5920597}{14408350}\bigg(a\,c\,p^{2} \bigg)^2
+\frac{1125167}{3159157}\bigg(c\,p \bigg)^2\\
> 0,
\end{dmath*}
\end{small}

\noindent where in each case, the second expression (i.e. the SOS decomposition) is obtained by an application of \texttt{`SumsOfSquares'} to the first expression. \hfill \fbox

\begin{rem}
To establish the inequality (\ref{ineq111m}), we need only consider Case 1. However, to show that the equality sign holds if and only if $X_1,X_2,X_3,X_4$ are independent, we must check all three cases.
\end{rem}

\section{Discussion}\setcounter{equation}{0}

The GPI conjecture is an extremely difficult problem to solve, particularly when some of the correlations are negative. The SOS method described in this paper can be used to rigorously verify \textit{any} specific case of the GPI, constrained only by computing power. Furthermore, as demonstrated in Section 4, this method is even powerful enough to prove GPIs with one exponent unbounded, a feat that is extremely difficult by purely theoretical methods (cf. \cite[Theorem 3.5]{RusSun}). On the other hand, should the GPI conjecture not hold in its full generality, our method may prove quite useful in the search for a counterexample. Our algorithm is efficient, straightforward and produces exact results. Furthermore, whereas calculations of multivariate Gaussian moments are often burdened by the constraints imposed by the covariance matrix, our method has the advantage of using free variables (with domain over the reals).

Despite the fact that the GPI is widely believed to be true, as of yet there has not been much to strongly support this presumption. Theorems \ref{thm34m} and \ref{thm111m} constitute never-before obtained results. In fact, we were able to verify many more GPIs with one exponent unbounded and higher fixed exponents, the proofs of which we omit seeing as the SOS decompositions, although exact, can be quite long. Thus, with the help of software, our work provides some of the first legitimate and substantial support to the correctness of the GPI. Therefore, we propose a stronger version of the GPI:
\begin{con}\label{conSOS} Let $n\ge 3$, $m_1,\dots,m_n\in\mathbb{N}$, and  $U_1,\dots,U_{n}$ be independent standard Gaussian random variables. Define the polynomial
$H$ on $\mathbb{R}^{n(n-1)/2}$ by
\begin{eqnarray*}
&&H(x_{21},x_{31},x_{32},\dots,x_{n1},\dots,
x_{n, n-1})\\
&=&E\left[U_1^{2m_1}\prod_{i=2}^{n}\left(\sum_{j=1}^{i-1}x_{ij}U_j+U_{i}\right)^{2m_i}\right]-(2m_1-1)!!\prod_{i=2}^{n}E\left[\left(\sum_{j=1}^{i-1}x_{ij}U_j+U_{i}\right)^{2m_i}\right].
\end{eqnarray*}
Then, $H$ has an SOS representation.
\end{con}

By approximation, we find that if Conjecture \ref{conSOS} is true then the GPI (\ref{14AA}) holds. By virtue of this conjecture, we have connected the probability inequality to  analysis, algebra, geometry, combinatorics, mathematical programming and computer science. In particular, research on SOS is abundant and ongoing (see, for example, \cite{LasBook}), but so far, no theoretical results in that field have been used to prove a GPI. Although software is used, our rigorous proofs establish the first meaningful link between the GPI and SOS, and should stimulate those working on SOS.

It is well-known that a non-negative multivariate polynomial  may not have an SOS representation. Denote by $P_{n,2d}$ the set of all non-negative polynomials in $n$ variables of degree at most $2d$ and $\sum_{n,2d}$ the set of all polynomials in $P_{n,2d}$ that are SOS. In 1888, Hilbert proved that $\sum_{n,2d}=P_{n,2d}$ if and only if $n=1$ or $d=1$ or $(n,d)=(2,2)$ (see \cite{Hilb}). However, we would like to point out that $E \left[\prod_{j=1}^{n}X_j^{2m_j}\right]$ itself is an SOS. To see this, let $U_1,\dots,U_n$ be  independent standard Gaussian random variables. For $x_{11},\dots,x_{1n},\dots,x_{n1},\dots,x_{nn}\in\mathbb{R}$, we can write:
$$
X_i=\sum_{j=1}^nx_{ij}U_j,\ \ \ \ 1\le i\le n.
$$
Then, we have
$$
\prod_{j=1}^{n}X_j^{2m_j}= a'b,
$$
where $a'=(1,x_{11},\dots,x_{nn},x_{11}^2,x_{11}x_{12},\dots,x_{1n}^{m_1}\cdots x_{nn}^{m_n})$ is a vector of monomials of length $l\leq n^{\sum_{j=1}^nm_j}$ and $b'=(P_1(U_1,\dots,U_n),P_2(U_1,\dots,U_n),\dots,P_l(U_1,\dots,U_n))$ is a vector of polynomials also of length $l$. Thus,
$$
E \left[\prod_{j=1}^{n}X_j^{2m_j}\right]=E \left[\left(\prod_{j=1}^{n}X_j^{m_j}\right)^2\right]=E[a'bb'a]=a'Qa,
$$
where $Q=bb'$ is a non-negative definite matrix. Therefore, $E \left[\prod_{j=1}^{n}X_j^{2m_j}\right]$ has an SOS representation.

\vskip 1cm

\begin{large} \noindent\textbf{Acknowledgements} \end{large} We thank Dr. Thomas Royen for fruitful discussion and encouragement with regards to the SOS method. We thank Dr. Victor Magron for pointing out Scheiderer's result \cite[Theorem 2.1]{Sch16} among other helpful suggestions. This work was supported by the Natural Sciences and Engineering Research Council of Canada (Nos. 559668-2021 and
4394-2018).

\end{document}